\documentclass[a4paper, 12pt]{amsart}

\usepackage[english]{babel}
\usepackage[latin1]{inputenc}
\usepackage[T1]{fontenc}
\usepackage{amsmath}
\usepackage{amsfonts}
\usepackage{amscd}

\theoremstyle{plain}
\newtheorem{theorem}{Theorem}[section]
\newtheorem{lemma}[theorem]{Lemma}
\newtheorem{proposition}[theorem]{Proposition}
\newtheorem{corollary}[theorem]{Corollary}
\newtheorem{question}[theorem]{Question}
\theoremstyle{definition}
\newtheorem{definition}[theorem]{Definition}

\newtheorem{example}[theorem]{Example}

\newcommand{\ZZ}{\mathbb{Z}} 
\newcommand{\CC}{\mathbb{C}} 
\newcommand{\st}{\;\vline\;} 

\newcommand{\tensor}{\otimes}

\newcommand{\iso}{\cong}      
\newcommand{\PP}{\mathbb{P}}  
\newcommand{\sheaf}[1]{\mathcal{#1}} 
\newcommand{\OO}{\sheaf{O}}   
\newcommand{\rto}{\dashrightarrow} 
\newcommand{\abdual}[1]{\widehat{#1}} 
\newcommand{\abddual}[1]{\ooalign{\raisebox{0.4ex}{$\widehat{\widehat{\ \ }}$}\cr$#1$}} 
\newcommand{\dual}[1]{{#1}^\vee}
\newcommand{\rtoiso}{\widetilde{\dashrightarrow}}
\newcommand{\toiso}{\widetilde{\longrightarrow}}

\DeclareMathOperator{\Pic}{Pic} 
\DeclareMathOperator{\NS}{NS} 

\author{Martin G. Gulbrandsen}
\address{Department of Mathematics\\P.O. Box 1053\\NO-0316
  Oslo\\Norway}
\email{martingu@math.uio.no}

\title{Lagrangian fibrations on generalized Kummer varieties}

\begin{document}

\begin{abstract}
  We investigate the existence of Lagrangian fibrations on the
  generalized Kummer varieties of Beauville. For a principally
  polarized abelian surface $A$ of Picard number one we find the
  following: The Kummer variety $K^nA$ is birationally equivalent to
  another irreducible symplectic variety admitting a Lagrangian
  fibration, if and only if $n$ is a perfect square. And this is the
  case if and only if $K^nA$ carries a divisor with vanishing
  Beauville-Bogomolov square.
\end{abstract}

\maketitle

\section{Introduction}

Let $X$ denote a projective irreducible symplectic variety of
dimension $2n$. We refer the reader to Huybrechts \cite{huybrechts}
for definitions and general background material. Matsushita
\cite{matsushita1, matsushita2} studied \emph{fibrations} of $X$, that
is, proper maps
\begin{equation}\label{eq:fib}
f\colon X\to B,
\end{equation}
such that a generic fibre is connected and has positive dimension.
Assuming $B$ to be projective and nonsingular, Matsushita showed that
$B$ has dimension $n$, its Hodge numbers agree with those of $\PP^n$,
and furthermore, a generic fibre of $f$ is a Lagrangian subvariety of $X$
of pure dimension $n$. It is a conjecture that $B$ is in fact
isomorphic to $\PP^n$.

The setup can be generalized slightly:

\begin{definition}\label{def:ratfib}
With $X$ and $B$ as above, a \emph{rational fibration} of $X$ over $B$
is a rational map
\begin{equation*}
f\colon X\rto B
\end{equation*}
such that there exists another projective irreducible symplectic
variety $X'$ and a birational equivalence $g\colon X'\rtoiso X$ such
that the composition $f\circ g$ is a (regular) fibration of $X'$
over $B$.
\end{definition}

A basic tool in the study of irreducible symplectic varieties is the
\emph{Beauville-Bogomolov form}, which is an integral quadratic form
$q$ on $H^2(X,\ZZ)$, satisfying
\begin{equation*}
q(\alpha)^n = c\deg(\alpha^{2n})
\end{equation*}
for some constant $c$. A birational map between
irreducible symplectic varieties induces an isomorphism on $H^2(-,\ZZ)$,
compatible with the Beauville-Bogomolov forms. It follows that in the
situation of definition \ref{def:ratfib}, the pullback $D=f^*H$ of any
divisor $H$ on $B$ satisfies $q(D) = 0$.
Conversely, one may ask:

\begin{question}\label{q:fib}
  Suppose $X$ carries a divisor $D$ with vanishing Beauville-Bogomolov
  square. Does $X$ admit a rational fibration over $\PP^n$?
\end{question}

One may try to answer the question for the known examples of
projective irreducible symplectic varieties. There are two standard
series of examples, both due to Beauville \cite{beauville}: The first
is the Hilbert scheme $S^{[n]}$ (of dimension $2n$) parametrizing
finite subschemes of length $n$ of a $K3$ surface $S$. The second is
the (generalized) \emph{Kummer variety} $K^nA$ (of dimension $2n-2$)
associated to an abelian surface $A$, defined as the fibre of
the map
\begin{equation}\label{eq:sum}
\sigma\colon A^{[n]} \to A.
\end{equation}
induced by the group law on $A$. The map $\sigma$ is locally trivial
in the étale topology, and in particular all fibers are isomorphic. So
there is no ambiguity in this definition. Recently, Sawon \cite{sawon}
and Markushevich \cite{markushevich} answered question \ref{q:fib} in
the affirmative for the Hilbert scheme $S^{[n]}$ of a generic K3
surface.  In this text, we consider the case of the Kummer varieties.
We obtain:

\begin{theorem}\label{thm:fib}
Let $A$ be an abelian surface carrying a curve $C\subset A$ of genus
$n+1$ with $n>2$. Then the Kummer variety $K^nA$ admits a rational
fibration
\begin{equation}\label{eq:ratfib}
f\colon K^nA\rto |\abdual{C}|\iso \PP^{n-1}
\end{equation}
where $\abdual{C}\in\Pic(\abdual{A})$ is a canonically defined dual
divisor class.
\end{theorem}

The divisor class $\abdual{C}$ is defined in example
\ref{ex:dualcurve}. The theorem is proved in section
\ref{sec:construction}. We have the following corollary, which answers
question \ref{q:fib} in the affirmative for the Kummer
varieties associated to a generic principally polarized abelian
surface, and which is proved in section \ref{sec:principal}:

\begin{corollary}\label{cor:principal}
  If the abelian surface $A$ has Picard number one and admits a
  principal polarization, then the following are equivalent, for each
  $n>2$:
\begin{enumerate}
\item The Kummer variety $K^nA$ admits a
rational fibration over $\PP^{n-1}$.
\item $K^nA$ carries a divisor with vanishing Beauville-Bogomolov
  square.
\item $n$ is a perfect square.
\end{enumerate}
\end{corollary}

Our work has been carried out independently of the works of Sawon
and Markushevich, but the construction is similar, except that
the Kummer variety case is easier, and doesn't require the use of
moduli spaces of twisted sheaves.

I would like to thank Geir Ellingsrud for numerous fruitful
discussions, and Manfred Lehn for introducing me to the question
of existence of Lagrangian fibrations.

\section{Preparation}

We work in the category of noetherian schemes over $\CC$. By a map of
schemes we mean a morphism in this category. By a sheaf on a scheme
$X$ we mean a coherent $\OO_X$-module.

If $A$ is an abelian variety, we denote the identity element for the
group law on $A$ by $0$, and if $a$ is a point on $A$, we write
$T_a\colon A\to A$ for translation by $a$. We write $\abdual{A}$ for
the dual abelian variety. We denote by $\sheaf{P}_x$ the homogeneous
line bundle on $A$ corresponding to a point $x\in \abdual{A}$. If $D$
is a divisor on $A$, we denote by
\begin{equation*}
\phi_D\colon A\to\abdual{A}
\end{equation*}
the map that takes a point $a\in A$ to the divisor $T_a^*D-D$.

In this section, we recall a few results from the literature.

\subsection{The Fourier-Mukai transform}

Let $X\to T$ be an abelian scheme over $T$,
and let $\abdual{X}\to T$ denote its dual abelian scheme. Let
$\sheaf{P}$ be the Poincaré line bundle on $X\times_T\abdual{X}$,
normalized such that the restrictions of $\sheaf{P}$ to $X\times 0$
and $0\times\abdual{X}$ are trivial. Let
\begin{equation*}
X \xleftarrow{p} X\times_T\abdual{X} \xrightarrow{q} \abdual{X}
\end{equation*}
denote the two projections.

Following Mukai \cite{mukai, mukai2} we define a functor $\abdual{S}$ from the
category of $\OO_X$-modules to the category of
$\OO_{\abdual{X}}$-modules by
\begin{equation*}
\abdual{S}(\sheaf{E}) = q_*(p^*(\sheaf{E})\tensor\sheaf{P}).
\end{equation*}
Reversing the roles of $X$ an $\abdual{X}$, we get a functor $S$
taking $\OO_{\abdual{X}}$-modules to $\OO_X$-modules.

\begin{definition}
An $\OO_X$-module $\sheaf{E}$ satisfies the \emph{weak index theorem} (WIT) with index $i$ if
\begin{equation*}
R^p\abdual{S}(\sheaf{E}) = 0\quad\text{for all $p\ne i$}.
\end{equation*}
The \emph{Fourier-Mukai transform} of such a sheaf $\sheaf{E}$ is the
$\OO_{\abdual{X}}$-module
\begin{equation*}
\abdual{\sheaf{E}} = R^i\abdual{S}(\sheaf{E}).
\end{equation*}
\end{definition}

For each $t\in T$, we may view $\sheaf{E}\tensor k(t)$ as a sheaf on
the fibre $X_t$, which is an abelian variety. We have the following
base change result:

\begin{theorem}[Mukai \cite{mukai2}]\label{thm:fm-basechange}
Let $\sheaf{E}$ be a sheaf on $X\to T$, flat over $T$.
The locus of points $t\in T$ such that $\sheaf{E}\tensor k(t)$
satisfies WIT is open. If $\sheaf{E}\tensor k(t)$ satisfies WIT with
index $i$ for all $t\in T$, then $\sheaf{E}$ also satisfies WIT with
index $i$, $\abdual{\sheaf{E}}$ is flat over $T$ and we have
\begin{equation*}
\abdual{\sheaf{E}\tensor k(t)}\iso\abdual{\sheaf{E}}\tensor k(t)
\end{equation*}
for all $t\in T$.
\end{theorem}

We will apply this only in the case $X=A\times T$, where $A$ is
an abelian surface, and view $\sheaf{E}$ as a family of sheaves on $A$
parametrized by $T$.

Mukai's discobery was the following:

\begin{theorem}[Mukai \cite{mukai}]\label{thm:fm-equiv}
Let $A$ be an abelian variety of dimension $g$. The functor
$\abdual{S}$ induces an equivalence of derived categories
\begin{equation*}
R\abdual{S}\colon D(A)\to D(\abdual{A})
\end{equation*}
with quasi-inverse taking a complex $\sheaf{K}^\bullet$ to
$(-1)_A^*RS(\sheaf{K}^\bullet)[g]$.

In particular, if $\sheaf{E}$ is a sheaf on $A$ satisfying WIT with
index $i$, then $\abdual{\sheaf{E}}$ also satisfies WIT, with index
$g-i$, and we have a natural isomorphism
\begin{equation*}
\abddual{\sheaf{E}}\iso (-1)_A^*\sheaf{E}.
\end{equation*}
\end{theorem}

A similar statement holds for an arbitrary abelian scheme
\cite{mukai2}, but we will not need this.

Mukai \cite{mukai2} also calculated the Chern character in
$H^*(\abdual{A},\ZZ)$ of $R\abdual{S}(\sheaf{K}^\bullet)$ in terms of the Chern
character of the complex $\sheaf{K}^\bullet$: There is a canonical
duality between the cohomology groups of $A$ and those of
$\abdual{A}$. Thus, using Poincaré duality, we may identify
\begin{equation}\label{eq:poincare}
H^p(\abdual{A},\ZZ)\iso H^{2g-p}(A,\ZZ).
\end{equation}
Writing $ch^p$ for the $2p$'th component of the Chern character, and
suppressing the isomorphism \eqref{eq:poincare}, Mukai found
\begin{equation}\label{eq:chern}
ch^p(R\abdual{S}(\sheaf{K}^\bullet)) = (-1)^p
ch^{g-p}(\sheaf{K}^\bullet).
\end{equation}
In particular, whenever $\sheaf{E}$ satisfies WIT with index $i$, we
have
\begin{equation}\label{eq:witchern}
ch^p(\abdual{\sheaf{E}}) = (-1)^{i+p}ch^{g-p}(\sheaf{E}).
\end{equation}
We remark that on an
abelian surface, the components of the Chern character are the rank,
the first Chern class and the Euler characteristic:
\begin{align*}
ch^0(\sheaf{E}) &= r(\sheaf{E}) & ch^1(\sheaf{E}) &= c_1(\sheaf{E}) & ch^2(\sheaf{E}) &= \chi(\sheaf{E})
\end{align*}

\begin{example}\label{ex:dualcurve}
Let $C\subset A$ be an effective curve of genus $g>1$ on an abelian
surface. By \cite[sec.~16]{mumford}, we have
\begin{equation*}
H^p(A,\sheaf{P}_x(C))=0,\quad \text{for all $p>0$ and all
  $x\in\abdual{A}$}.
\end{equation*}
Hence, $\OO_A(C)$ satisfies WIT with index $0$, and
$\abdual{\OO_A(C)}$ is locally free on $\abdual{A}$. Applying
formula \eqref{eq:witchern}, we see that
$\abdual{\OO_A(C)}$ has rank equal to $\chi(\OO_A(C)) = g-1$ and
Euler characteristic $1$, whereas under the isomorphism
\eqref{eq:poincare}, we have
\begin{equation*}
c_1(\abdual{\OO_A(C)}) = - [C].
\end{equation*}
Thus, defining $\abdual{C}\in\Pic(\abdual{A})$ to be the divisor
class such that
\begin{equation*}
\OO_A(-\abdual{C}) \iso \det \abdual{\OO_A(C)},
\end{equation*}
we see that $[C]$ corresponds to $[\abdual{C}]$ under the isomorphism
\eqref{eq:poincare}.
\end{example}

\subsection{Moduli of sheaves on an abelian surface}\label{sec:moduli}

Let $A$ be an abelian surface and fix a polarization $H$. By a (semi-)
stable sheaf on $A$ we will mean a Gieseker (semi-) stable sheaf with
respect to $H$.  Fixing a rank $r\ge 0$, first Chern class
$c_1\in\NS(A)$ and Euler characteristic $\chi$, we denote by
$M_A(r,c_1,\chi)$ the Simpson moduli space of stable sheaves with
the given invariants. In the cases of interest to us, stability and
semi-stability will be equivalent, so that $M_A(r,c_1,\chi)$ is going
to be projective.

We will in fact only consider sheaves of rank one or zero, so
we note that every torsion free sheaf of rank one is stable, whereas
in the rank zero case, Riemann-Roch gives the following condition: A
pure one dimensional sheaf $\sheaf{E}$ on $A$ is stable if and only if
we have
\begin{equation}\label{eq:torsionstab}
\frac{\chi(\sheaf{F})}{\deg_H(\sheaf{F})} <
\frac{\chi(\sheaf{E})}{\deg_H(\sheaf{E})}
\end{equation}
for every non trivial proper subsheaf
$\sheaf{F}\subset\sheaf{E}$.

Yoshioka \cite{yoshioka} defines a (regular) map
\begin{equation}\label{eq:albanese}
\alpha\colon M_A(r,c_1,\chi)\to A\times\abdual{A}
\end{equation}
that can be described at the level of sets as follows, except that we
take the liberty to make a sign change: Choose a representative
$\sheaf{L}\in\Pic(A)$ in the class $c_1$, and also a representative
$\sheaf{L}'\in\Pic(\abdual{A})$ in the class corresponding to $c_1$
via Poincaré duality \eqref{eq:poincare}. Then define
$\alpha=(\delta,\abdual{\delta})$, where
\begin{align*}
\delta(\sheaf{F}) &= \det(R\abdual{S}(\sheaf{F}))^{-1}\tensor\sheaf{L}'^{-1}\\
\abdual{\delta}(\sheaf{F}) &= \det(\sheaf{F})\tensor\sheaf{L}^{-1}.
\end{align*}
Note that $\abdual{\delta}(\sheaf{F})$ is an element of
$\Pic^0(A)=\abdual{A}$ and, by \eqref{eq:chern},
$\delta(\sheaf{F})$ is an element of $\Pic^0(\abdual{A})=A$.

Whenever $\dim M_A(r,c_1,\chi) \ge 8$, Yoshioka obtains the following:

\begin{theorem}[Yoshioka \cite{yoshioka}]\label{thm:yoshioka}
Assume the polarization $H$ is generic.
\begin{enumerate}
\item $M_A(r,c_1,\chi)$ is
  deformation equivalent to $A^{[n]}\times\abdual{A}$ for suitable
  $n$.
\item The map $\alpha$ in \eqref{eq:albanese} is locally trivial in the
étale topology.
\item A fibre $K_A(r,c_1,\chi)$ of the map $\alpha$ is deformation
  equivalent to the Kummer variety $K^nA$.
In particular, $K_A(r,c_1,\chi)$ is an irreducible symplectic
variety. 
\end{enumerate}
\end{theorem}

As we will be free to choose the polarization $H$ arbitrarily, the
genericity hypothesis will not be of importance to us. We remark,
however, that in the case where $A$ has Picard number one, every
polarization is generic.

\subsection{The Beauville-Bogomolov form on Kummer varieties}

The Beauville-Bogomolov form can be explicitly described in the case
of Kummer varieties. The description of $H^2(K^nA,\CC)$ is due to
Beauville \cite{beauville}, and the calculation of the
Beauville-Bogomolov form can be found in \cite[sec.~4.3.1]{yoshioka} or \cite[prop.~1]{britze}.

Firstly, there is a canonical monomorphism
\begin{equation}\label{eq:H2mono}
H^2(A,\CC)\to H^2(K^nA,\CC)
\end{equation}
which is compatible with the Hodge
structure. Secondly, there is a primitive integral class $\epsilon\in
H^2(K^nA,\CC)$ such that $2\epsilon$ is the fundamental class of the
locus $E\subset K^nA$ consisting of non reduced subschemes. Thus
$\epsilon$ is a  $(1,1)$-class. Together, $H^2(A,\CC)$ and $\epsilon$
generate $H^2(K^nA,\CC)$. In fact, we have:

\begin{proposition}
There is a direct sum decomposition
\begin{equation*}
H^2(K^nA,\CC) \iso H^2(A,\CC)\oplus\CC\epsilon
\end{equation*}
which is orthogonal with respect to the Beauville-Bogomolov form
$q$. Furthermore, the restriction of $q$ to $H^2(A,\CC)$ is the
intersection form on $A$, whereas
\begin{equation*}
q(\epsilon) = -2n.
\end{equation*}
\end{proposition}

We are interested in classes in $H^2(K^nA,\CC)$ coming from
divisors, that is, the Neron-Severi group $\NS(K^nA)$.
Since the inclusion \eqref{eq:H2mono} is compatible with the Hodge
structure, and $\epsilon$ is a primitive $(1,1)$-class, we find
\begin{equation}\label{eq:NS}
NS(K^nA)\iso NS(A)\oplus\ZZ\epsilon,
\end{equation}
by the Lefschetz theorem on $(1,1)$-classes.

\section{Construction}\label{sec:construction}

Consider the setup of theorem \ref{thm:fib}, that is, we have a
curve $C\subset A$ of genus $n+1$ on an abelian surface $A$.

To construct the fibration in theorem \ref{thm:fib}, we
want to associate to each $\xi\in A^{[n]}$ a curve in a certain linear
system. As a first try, one might ask whether there exists a curve in
the linear system $|C|$ containing $\xi$. This turns out to be too
restrictive:
\begin{lemma}\label{lem:nonincidence}
A generic element $\xi\in A^{[n]}$ is not contained in any curve in the
linear system $|C|$.
\end{lemma}

\begin{proof}
As we have seen in example \ref{ex:dualcurve}, we have
\begin{equation*}
H^p(A,\OO_A(C)) = 0\quad \text{for all $p>0$}
\end{equation*}
and thus, by Riemann Roch,
\begin{equation*}
\dim H^0(A,\OO_A(C)) = \chi(\OO_A(C)) = n.
\end{equation*}
Thus the complete linear system $|C|$ has dimension $n-1$. It follows
that the set of subschemes $\xi\in A^{[n]}$ contained in a curve in
$|C|$ forms a family of dimension $2n-1$. On the other hand, $A^{[n]}$
has dimension $2n$.
\end{proof}

Let us, starting from the observation in the lemma, sketch our
construction: By allowing not only curves in $|C|$, but in the linear
systems associated to $\sheaf{P}_x(C)$ for any $x\in\abdual{A}$, we
see that we ``win'' two more degrees of freedom: The set of length $n$
subschemes contained in a curve in $|\sheaf{P}_x(C)|$, for some
$x\in\abdual{A}$, forms a family
of dimension $2n+1$.  Since, again, $A^{[n]}$ has dimension $2n$, we
expect the locus
\begin{equation}\label{eq:D}
D_\xi = \{x\in\abdual{A}\st
H^0(A,\sheaf{I}_\xi\tensor\sheaf{P}_x(C))\ne 0\}
\end{equation}
to be a curve. We will see that this is indeed true for generic
$\xi$, and furthermore, when $\xi$ is a generic element of the Kummer
variety $K^nA$, the curve $D_\xi$ belongs to the linear system
$|\abdual{C}|$. The fibration $f$ in theorem \ref{thm:fib} is
given by sending $\xi$ to $D_\xi$.

More precisely we will see that, for generic $\xi\in K^nA$, the sheaf
$\sheaf{I}_\xi(C)$ satisfies WIT with index $1$. Sending $\xi$ to the
Fourier-Mukai transform $\abdual{\sheaf{I}_\xi(C)}$ induces a
birational equivalence
\begin{equation}\label{eq:fib1}
K^nA \rtoiso K_{\abdual{A}}(0,[\abdual{C}],-1)
\end{equation}
where the target space is the symplectic variety introduced in
Yoshioka's theorem \ref{thm:yoshioka}. The sheaves
parametrized by $K_{\abdual{A}}(0,[\abdual{C}],-1)$ are supported on
curves in the linear system $|\abdual{C}|$, and sending a sheaf to its
support defines a map
\begin{equation}\label{eq:fib2}
K_{\abdual{A}}(0,[\abdual{C}],-1) \to |\abdual{C}|.
\end{equation}
The composition of the two maps \eqref{eq:fib1} and \eqref{eq:fib2}
again gives us the fibration of theorem \ref{thm:fib}. We remark that
the support of $\abdual{\sheaf{I}_\xi(C)}$ is precisely the curve
$D_\xi$ in \eqref{eq:D}. In fact, the fibers of
$\abdual{\sheaf{I}_\xi(C)}$ are the vector spaces
\begin{equation*}
\abdual{\sheaf{I}_\xi(C)}\tensor k(x) \iso H^1(A,\sheaf{I}_\xi\tensor\sheaf{P}_x(C))
\end{equation*}
which vanish precisely when
$H^0(A,\sheaf{I}_\xi\tensor\sheaf{P}_x(C))$ vanish, since both the
Euler characteristic and the second cohomology of
$\sheaf{I}_\xi\tensor\sheaf{P}_x(C)$ is zero.

It turns out to be convenient to extend the setup as follows: We will
first see that there is a natural identification
$A^{[n]}\times\abdual{A}\iso M_A(1,[C],0)$ in such a way that the Kummer
variety is recovered as the fibers of the map
\begin{equation*}
\alpha\colon M_A(1,[C],0)\to A\times\abdual{A}
\end{equation*}
introduced in section \ref{sec:moduli}. Then we will construct a commutative diagram
\begin{equation}\label{eq:diagram}
\begin{CD}
M_A(1,[C],0)      & \stackrel{G}{\rtoiso}   & M_{\abdual{A}}(0,[\abdual{C}],-1) @>>> P    \\
@VV{\alpha}V                                  @VV{\alpha}V                           @VVV \\
A\times\abdual{A} & \stackrel{\psi}{\toiso} & \abdual{A}\times A                @>>> \abdual{A}
\end{CD}
\end{equation}
where $G$ is a birational map induced by the Fourier-Mukai transform,
$\psi$ is an isomorphism and $P\to\abdual{A}$ is a projective bundle
with the complete linear system associated to $\sheaf{P}_x(C)$ as the
fibre over $x$. Restricting the whole diagram to a fibre, we
recover the maps \eqref{eq:fib1} and \eqref{eq:fib2}.

\subsection{Rank one sheaves and the Hilbert scheme}

As usual, $A^{[n]}$ can be regarded as a moduli space of rank one
sheaves on $A$. More precisely, there is an isomorphism
\begin{equation}\label{eq:hilbmod-iso}
A^{[n]}\times\abdual{A}\iso M_A(1,0,-n)
\end{equation}
which, on the level of sets, is given by the map
\begin{equation*}
(\xi, x)\mapsto \sheaf{I}_\xi\tensor\sheaf{P}_x.
\end{equation*}

By twisting with $C$, we can furthermore identify $M_A(1,0,-n)$ with
$M_A(1,[C],0)$. Including the isomorphism \eqref{eq:hilbmod-iso}, we
can thus identify
\begin{equation*}
A^{[n]}\times\abdual{A} \iso M_A(1,[C],0).
\end{equation*}
We want to describe the composition
\begin{equation*}
A^{[n]}\times\abdual{A}\iso M_A(1,[C],0)\xrightarrow{\alpha}
A\times\abdual{A}
\end{equation*}
where $\alpha$ is the map \eqref{eq:albanese} of Yoshioka.  Recall
that to define $\alpha$, we must choose invertible sheaves $\sheaf{L}$
and $\sheaf{L}'$ representing $c_1=[C]$ on $A$ and on $\abdual{A}$,
respectively. By example \ref{ex:dualcurve}, we have the natural choices
\begin{align*}
\sheaf{L} &= \OO_A(C) & \sheaf{L}' &= \OO_{\abdual{A}}(\abdual{C}),
\end{align*}
and then we have:

\begin{lemma}
The diagram
\begin{equation*}
\begin{CD}
A^{[n]}\times\abdual{A}           & \iso       & M_A(1,[C],0)\\
@VV{\sigma\times 1_{\abdual{A}}}V              @VV{\alpha}V\\
A\times\abdual{A}                 @>{\theta}>> A\times\abdual{A}
\end{CD}
\end{equation*}
is commutative, where $\theta$ is the
isogeny
\begin{equation*}
\theta(a,x) = (a+\phi_{\abdual{C}}(x),x).
\end{equation*}
In particular, the fibers $K^nA$ on the left are
taken isomorphically to the fibers $K_A(1,[C],0)$ on the right.
\end{lemma}

\begin{proof}
  Let us, for the sake of readability, use additive notation in the
  Picard groups. Firstly, we have
\begin{equation*}
\abdual{\delta}(\sheaf{I}_\xi\tensor\sheaf{P}_x(C)) =
\det(\sheaf{I}_\xi\tensor\sheaf{P}_x(C)) + \OO_A(-C) = \sheaf{P}_x.
\end{equation*}
Secondly, applying the Fourier-Mukai functor to the short exact sequence
\begin{equation}\label{eq:xi-ses}
0\to\sheaf{I}_\xi\tensor\sheaf{P}_x(C)\to\sheaf{P}_x(C)\to \OO_\xi\to 0
\end{equation}
we obtain an exact sequence
\begin{equation*}
0\to \abdual{S}(\sheaf{I}_\xi\tensor\sheaf{P}_x(C))
 \to \abdual{S}(\sheaf{P}_x(C))
 \to \abdual{S}(\OO_\xi)
 \to R^1\abdual{S}(\sheaf{I}_\xi\tensor\sheaf{P}_x(C)) \to 0,
\end{equation*}
since $\sheaf{P}_x(C)$ satisfies WIT with index $0$, as in example
\ref{ex:dualcurve}. Thus we have
\begin{equation*}
\delta(\sheaf{I}_\xi\tensor\sheaf{P}_x(C))
= - \det \abdual{S}(\sheaf{P}_x(C)) + \det \abdual{S}(\OO_\xi) + \OO_{\abdual{A}}(-\abdual{C}).
\end{equation*}
By direct computation, we find that the Fourier-Mukai transform
$\abdual{S}(\OO_\xi)$ is the direct sum
$\bigoplus_{a\in\xi}\sheaf{P}_a$, where the points $a\in\xi$ should be
repeated according to their multiplicity. Hence
\begin{equation*}
\det(\abdual{S}(\OO_\xi)) =
\sheaf{P}_{\sigma(\xi)}
\end{equation*}
where $\sigma$ is the summation map \eqref{eq:sum}.
Furthermore, by \cite[sec.~3.1]{mukai}, tensoring with $\sheaf{P}_x$ before
applying $\abdual{S}$ is the same thing as translating with $x$ after
applying $\abdual{S}$, and hence
\begin{equation*}
\det(\abdual{S}(\sheaf{P}_x(C)) = \OO_A(-T_x^*\abdual{C})
\end{equation*}
by the definition of $\abdual{C}$ in example
\ref{ex:dualcurve}. 
Thus
\begin{equation*}
\delta(\sheaf{I}_\xi\tensor\sheaf{P}_{\sigma(\xi)}(C)) =
\sheaf{P}_{\sigma(\xi)}+\OO_A(T_x^*\abdual{C} - \abdual{C}).
\end{equation*}
More concisely, we may write this as
\begin{equation*}
\alpha(\sheaf{I}_\xi\tensor\sheaf{P}_x(C)) = (\sigma(\xi) + \phi_{\abdual{C}}(x),x).
\end{equation*}
which is what we wanted to prove.
\end{proof}

\subsection{The weak index property}

\begin{lemma}\label{lem:wit}
The (open) locus of sheaves $\sheaf{E}\in M_A(1,[C],0)$ satisfying WIT
with index $1$ is non empty. In fact, there exists WIT sheaves in every
fibre $K_A(1,[C],0)$ of $\alpha$.
\end{lemma}

\begin{proof}
By \cite[sec.~3.1]{mukai}, the operations of translation and twisting by
a homogeneous line bundle
\begin{equation*}
\sheaf{E}\mapsto
T_a^*\sheaf{E},\quad\sheaf{E}\mapsto\sheaf{E}\tensor\sheaf{P}_x
\end{equation*}
do not affect the WIT'ness of a sheaf $\sheaf{E}$. Thus, it is enough
to prove the existence of a WIT sheaf in $M_A(1,[C],0)$, since we can move
such a sheaf to any fibre of $\alpha$ by translating and twisting.

Let $\sheaf{E}=\sheaf{I}_\xi(C)$. We have
\begin{equation*}
H^2(\abdual{A},\sheaf{I}_\xi(C)\tensor\sheaf{P}_x) = 0
\end{equation*}
for all $x\in\abdual{A}$, for instance by the short exact sequence
\eqref{eq:xi-ses}, so $R^2\abdual{S}(\sheaf{I}_{\xi}(C))=0$. 

Furthermore, by lemma
\ref{lem:nonincidence}, we have
\begin{equation*}
H^0(A, \sheaf{I}_{\xi}(C)) = 0
\end{equation*}
for generic $\xi$. 
But $\abdual{S}(\sheaf{I}_{\xi}(C))$ is torsion free, hence we conclude
that
$\abdual{S}(\sheaf{I}_{\xi}(C))=0$ for generic $\xi$. Thus
$\sheaf{I}_\xi(C)$ satisfies WIT with index $1$.
\end{proof}

\subsection{Stability}

\begin{lemma}\label{lem:stable}
Let $\sheaf{E}$ be a sheaf in $M(1,[C],0)$ satisfying WIT with index
$1$. Then the Fourier-Mukai transform $\abdual{\sheaf{E}}$ is stable with
respect to any polarization of $\abdual{A}$.
\end{lemma}

\begin{proof}
We first show that $\abdual{\sheaf{E}}$ is pure. Being the
Fourier-Mukai transform of a WIT sheaf with index $1$,
$\abdual{\sheaf{E}}$ itself satisfies WIT with index $1$. It has rank
zero and first Chern class $[\abdual{C}]\ne 0$, hence it is one
dimensional. If $\sheaf{T}\subset\abdual{\sheaf{E}}$ is a zero
dimensional subsheaf, then $\sheaf{T}$ satisfies WIT with index $0$, but
\begin{equation*}
S(\sheaf{T}) \subseteq S(\abdual{\sheaf{E}}) = 0
\end{equation*}
and hence $\sheaf{T} = 0$. Thus $\abdual{\sheaf{E}}$ is pure of
dimension $1$.

Suppose $\sheaf{F}\subset\abdual{\sheaf{E}}$ were a destabilizing
subsheaf. Then $\sheaf{F}$ also satisfies WIT with index $1$. Since its
support, defined by the Fitting ideal, is contained in the support
of $\abdual{\sheaf{E}}$, and the degree of a one dimensional sheaf is given by
intersecting its (Fitting) support with the polarization, we have
\begin{equation*}
\deg(\sheaf{F})\le\deg(\abdual{\sheaf{E}})
\end{equation*}
with respect to any polarization of $\abdual{A}$. On the other hand,
since $\sheaf{F}$ is destabilizing, we have by \eqref{eq:torsionstab}
\begin{equation*}
\frac{\chi(\sheaf{F})}{\deg(\sheaf{F})} >
\frac{\chi(\abdual{\sheaf{E})}}{\deg(\abdual{\sheaf{E}})}
\end{equation*}
and thus
\begin{equation*}
\chi(\sheaf{F}) > \chi(\abdual{\sheaf{E}}) = -1.
\end{equation*}
Since the Fourier-Mukai transform $\abdual{\sheaf{F}}$ has rank
$-\chi(\sheaf{F}) < 1$ by equation \eqref{eq:witchern}, it must be a
torsion sheaf. Now, applying the Fourier-Mukai functor to the exact
sequence
\begin{equation*}
0\to\sheaf{F}\to\abdual{\sheaf{E}}\to\abdual{\sheaf{E}}/\sheaf{F}\to 0
\end{equation*}
we obtain a left exact sequence
\begin{equation*}
0\to S(\abdual{\sheaf{E}}/\sheaf{F})
 \to \abdual{\sheaf{F}}
 \to \abddual{\sheaf{E}} \iso (-1)^*\sheaf{E}
\end{equation*}
where theorem \ref{thm:fm-equiv} is applied to obtain the isomorphism
on the right. But both $S(\abdual{\sheaf{E}}/\sheaf{F})$ and
$(-1)^*\sheaf{E}$ are torsion free, hence it is impossible for the
middle term $\abdual{\sheaf{F}}$ to be torsion. Thus
we have reached a contradiction.
\end{proof}

We are now ready to construct the leftmost square in diagram
\eqref{eq:diagram}: Let $U\subset M_A(1,[C],0)$ denote the set of
sheaves satisfying WIT with index $1$. By theorem
\ref{thm:fm-basechange}, $U$ is open, and by lemma \ref{lem:wit}, $U$
is non empty. Let $\sheaf{U}$ denote the restriction of the universal
family on $M_A(1,[C],0)$ to $U$. Applying theorem
\ref{thm:fm-basechange} again, $\sheaf{U}$ satisfies WIT with index
$1$, and its Fourier-Mukai transform $\abdual{\sheaf{U}}$ is a flat
family of sheaves on $\abdual{A}$ parametrized by $U$. The fibers of
$\abdual{\sheaf{U}}$ are stable by lemma \ref{lem:stable}, and by
equation \eqref{eq:witchern} they have rank one, first Chern class
$[\abdual{C}]$ and Euler characteristic $-1$. Thus there is an induced
rational map
\begin{equation*}
G\colon M_A(1,[C],0) \rto M_{\abdual{A}}(0,[\abdual{C}],-1)
\end{equation*}
which is regular on $U$. In fact, by theorem \ref{thm:fm-equiv}, the
restricion of $G$ to $U$ is an isomorphism. As
$M_{\abdual{A}}(0,[\abdual{C}],-1)$ is irreducible by theorem
\ref{thm:yoshioka}, $G$ is birational. Let us verify that $G$ fits
into the diagram \eqref{eq:diagram}, i.e.~we check the commutativity
of the leftmost square. So let $\sheaf{E}$ be a sheaf in
$M_A(1,[C],0)$ satisfying WIT with index $1$. Then
\begin{align*}
\delta(\sheaf{E}) &= \det(\abdual{\sheaf{E}})\tensor\OO_{\abdual{A}}(-\abdual{C})\\
\abdual{\delta}(\sheaf{E}) &= \det(\sheaf{E})\tensor\OO_A(-C)
\end{align*}
whereas
\begin{align*}
\delta(\abdual{\sheaf{E}}) &=
\det(\abddual{\sheaf{E}})\tensor\OO_A(-C) =
(-1)_A^*\det(\sheaf{E})\tensor\OO_A(-C)\\
\abdual{\delta}(\abdual{\sheaf{E}}) &=
\det(\abdual{\sheaf{E}})\tensor\OO_{\abdual{A}}(-\abdual{C}).
\end{align*}
Thus we see that, defining the map $\psi$ in diagram
\eqref{eq:diagram} by
\begin{equation*}
\psi(a,x) = (-x,a) + ((-1)^*C-C, 0),
\end{equation*}
the left square in that diagram commutes. Since, by lemma
\ref{lem:wit}, no fibre $K_A(1,[C],0)$ of $\alpha$ is contained in the
base locus of $G$, we conclude that $G$ restricts to a birational
equivalence
\begin{equation*}
g\colon K_A(1,[C],0)\rtoiso K_{\abdual{A}}(0,[\abdual{C}],-1).
\end{equation*}

\subsection{The fibration}

Continuing example \ref{ex:dualcurve}, note that, by the base change
theorem, the fibre of $\abdual{\OO_A(C)}$ over $x\in\abdual{A}$ is
canonically ismorphic to $H^0(A, \sheaf{P}_x(C))$. Thus, the
associated projective bundle
\begin{equation}\label{eq:P}
P = \PP\left(\dual{\abdual{\OO_A(C)}}\right) \to \abdual{A}
\end{equation}
has the complete linear systems associated to $\sheaf{P}_x(C)$ as
fibers.

The Fitting ideal of a sheaf $\sheaf{F}$ in
$M_{\abdual{A}}(0,[\abdual{C}],-1)$ defines a curve representing the
first Chern class of $\sheaf{F}$, and hence a point in the bundle $P$.
The map of sets
\begin{equation}\label{eq:fitting}
F\colon M_{\abdual{A}}(0,[\abdual{C}],-1) \to P
\end{equation}
thus obtained is in fact a (regular) map of varieties, since formation
of the Fitting ideal commutes with base change. Clearly, $F$ fits into
diagram \eqref{eq:diagram}, making its rightmost square commute. Thus,
restricting $F$ to the fibre $K_{\abdual{A}}(0,[\abdual{C}],-1)$ above
zero in $A\times\abdual{A}$, we find a map
\begin{equation}
f\colon K_{\abdual{A}}(0,[\abdual{C}],-1)\to |\abdual{C}|.
\end{equation}

We claim that $f$ is a fibration, i.e.~a generic fibre is connected.
For this, let $D\in |\abdual{C}|$ be a nonsingular curve. Viewing $D$
as a point in $P$, the fibre $F^{-1}(D)$ is just the Jacobian
$J^{n-1}$ of $D$, parametrizing invertible sheaves of degree $n-1$ on
$D$. The restriction of $\alpha$ to $J^{n-1}$ can be identified with
the summation map
\begin{equation}\label{eq:jacsum}
J^{n-1} \to \abdual{A}
\end{equation}
sending a divisor $\sum n_i p_i$ on $D$ to the point $\sum n_i p_i$ on
$A$, using the group law on $\abdual{A}$. Thus, the fibre of $f$ above
$D$ equals a fibre of the map \eqref{eq:jacsum}, which is connected.
This concludes the proof of theorem \ref{thm:fib}.

\subsection{Principally polarized surfaces}\label{sec:principal}

Let $(A,H)$ satisfy the conditions of corollary \ref{cor:principal}.
The implication (1) $\implies$ (2) is automatic, as explained in the
introduction. For the implication (2) $\implies$ (3), suppose $K^nA$
admits a divisor $D$ with vanishing Beauville-Bogomolov square,
corresponding to $rH+s\epsilon$ under the isomorphism \eqref{eq:NS}.
Then
\begin{equation*}
0 = q(D) = (rH)^2 + s^2 q(\epsilon) = 2r^2 - 2s^2n
\end{equation*}
from which it is immediate that $n$ is a perfect square.

Finally, the implication (3) $\implies$ (1) follows from theorem
\ref{thm:fib}: If $n=m^2$ is a perfect square, the effective curve
$C=mH$ has genus $n+1$, and hence the theorem applies.
Corollary \ref{cor:principal} is proved.

\section{On the base locus}

We still consider a principally polarized abelian surface $(A,H)$ and
a perfect square $n=m^2$. Then there does exist $\xi\in A^{[n]}$ not
satisfying WIT: One can check that this is the case whenever $\xi\in
A^{[n]}$ is contained in some translate $T_a^{-1}(H)$ of the
polarization.  Preliminary investigation suggests that for $n=4$, the
birational equivalence \eqref{eq:fib1} induced by the Fourier-Mukai
transform is a Mukai elementary transformation, which cannot be
extended to an isomorphism. We will return to this question in a later
paper.

\bibliographystyle{amsplain}
\bibliography{kummerfib.bib}

\end{document}